\documentclass{amsart}
\usepackage{epsfig, amssymb}
\address{Department of Mathematics,
University of Toronto, Toronto ON M5S 2E4 Canada}

\newtheorem{theorem}{Theorem}[section]
\newtheorem{thm}{Theorem}

\newtheorem{proposition}[theorem]{Proposition}

\newtheorem{property}[theorem]{Property}

\theoremstyle{definition}
\newtheorem{definition}[theorem]{Definition}
\newtheorem{remark}[theorem]{Remark}
\newtheorem{example}[theorem]{Example}

\newcommand{\ignore}[1]{\relax}

\newcommand{\C}{\mathbb C}
\newcommand{\R}{\mathbb R}
\newcommand{\Z}{\mathbb Z}
\newcommand{\N}{\mathbb N}

\newcommand{\T}{\mathbb T}
\newcommand{\tp}{\mathbb T \mathbb P}
\newcommand{\U}{\mathcal U}

\newcommand{\M}{\mathcal M}

\newcommand{\val}{\operatorname{val}}

\newcommand{\pic}{\operatorname{Pic}}

\newcommand{\ver}{\operatorname{Vert}}

\newcommand{\vdim}{\operatorname{vdim}}

\newcommand{\tor}{(\C^\times)^{n}}

\newcommand{\ft}{\operatorname{ft}}
\newcommand{\ev}{\operatorname{ev}}

\newcommand{\am}{\mathcal{A}}

\newcommand{\cp}{{\mathbb C}{\mathbb P}}
\newcommand{\rp}{{\mathbb R}{\mathbb P}}

\newcommand{\Log}{\operatorname{Log}}

\renewcommand{\setminus}{\smallsetminus}

\newcommand{\Jac}{\operatorname{Jac}}
\newcommand{\Div}{\operatorname{Div}}
\newcommand{\Pic}{\operatorname{Pic}}
\newcommand{\Sym}{\operatorname{Sym}}

\newcommand{\OO}{\mathcal O}

\title{Tropical geometry and its applications}

\author[Grigory Mikhalkin]{Grigory Mikhalkin}

\begin{document}

\begin{abstract}
From a formal perspective tropical geometry can be viewed as a
branch of geometry manipulating with certain piecewise-linear
objects that take over the r\^ole of classical algebraic varieties.
This talk outlines some basic notions of this area and surveys some
of its applications for the problems in classical (real and complex)
geometry.
\end{abstract}

%\begin{classification}
%14A99, 14H50, 14N10, 52B20.
%\end{classification}
%
%\begin{keywords}
%Tropical geometry, amoebas, patchworking,
%enumerative geometry.
%\end{keywords}

\thanks{The author is grateful to
the Institut Henri Poincar\'e (Paris) and the IHES for hospitality during the preparation
of this talk.
Author's research is supported in part by NSERC.}
\maketitle

\section{Introduction}

From a geometric point of view tropical geometry describes worst
possible degenerations of the complex structure on an $n$-fold $X$.
Such degenerations cause $X$ to collapse onto an $n$-dimensional
(over $\R$, i.e. $(\frac{\dim X}{2})$-dimensional) base $B$ which is
a piecewise-linear polyhedral complex, see e.g \cite{KS} for a
conjectural picture in the special case of Calabi-Yau $n$-folds.

According to an idea of Kontsevich such degenerations can
be useful, in particular, for computations of the Gromov-Witten
invariants of $X$ as holomorphic curves degenerate to graphs $\Gamma\subset B$.
A similar picture appeared in the work of Fukaya (see e.g. \cite{Fu})
where graphs come as degenerations of holomorphic membranes.

Tropical geometry formalizes the base $B$ as an ambient variety
so that the graphs $\Gamma$ become curves in $B$.
Some problems in complex and real geometry then may be reduced
to problems of tropical geometry which are often much easier
to solve, thanks to the piecewise-linear nature of the subject.
Local considerations in tropical geometry correspond to some
standard models in classical geometry while the combinatorial
tropical structure encodes the way to glue these models together.
In this sense it may be viewed an extension of the Viro patchworking \cite{Vi}.

This talk takes a geometric point of view on tropical geometry and
surveys its basic notions
as well as some of its applications for classical algebro-geometrical problems.
The author is indebted to Ya. Eliashberg, M. Kontsevich,
A. Losev, B. Sturmfels and O. Viro
for many useful conversations on tropical geometry.
%that influenced the presentation of this talk.

\section{Tropical algebra}

\subsection{Tropical vs. classical arithmetics}

The term {\em tropical semirings} was reputedly
invented by a group of French computer scientists
to commemorate their Brazilian colleague Imre Simon.
For our purposes we use just one of these semirings, the
tropical semifield
$\T=\R\cup\{-\infty\}$
equipped with the operations of addition and multiplication
defined by
$$``a+b"=\max\{a,b\},\ ``ab"=a+b.$$
We use the quotation marks to distinguish the tropical operations
from the classical addition and multiplication on $\R\cup\{-\infty\}$.
It is easy to check that the tropical operations are commutative,
associative and satisfy to the distribution law.

There is no tropical subtraction as the operation $``+"$ is idempotent,
but we have the tropical division $``\frac{a}{b}"=a-b$, $b\neq-\infty$.
This makes $\T$ an idempotent semifield. Its additive zero is $0_\T=-\infty$,
its multiplicative unit is $1_\T=0$.

\begin{remark}
A classical example of a semifield is the semifield $\R_{\ge 0}$ of
nonnegative numbers with classical addition and multiplication.
We have the map $\log_t:\R_{\ge 0}\to \T$
%can be used to push forward the
%arithmetic operations from $\R_{\ge 0}$
between semifields $\R_{\ge0}$ for any $t>1$.
The larger is $t$ the closer is the map $\log_t$ to a homomorphism.
More specifically, $\log_t$ is an isomorphism up to the error of $\log_t(2)$
while $\lim\limits_{t\to+\infty}\log_t2=0$.
Indeed, we have $\log_t(ab)=``\log_t(a)\log_t(b)"$
and
$$``\log_t(a)+\log_t(b)"\le \log_t(a+b)\le ``\log_t(a)+\log_t(b)" + \log_t2$$
by elementary considerations.
Thus, $\T$ can be considered as the $t\to+\infty$ limit semifield
of a family of (mutually isomorphic) semifields $\R\cup\{-\infty\}$ equipped
with arithmetic operations induced from $\R_{\ge0}$ by $\log_t$.
In this sense $\T$ is considered to be the result of {\em dequantization}
of classical arithmetics by Maslov et al., see e.g. \cite{LiMa} or \cite{Vi-deq}
for a geometric version of this dequantization.
\end{remark}

\begin{remark}
The semifield $\T$ is closely related to non-Archimedean fields $K$ with (real)
valuations $\val:K\to\T.$
Indeed, $\val$ is a valuation if for any $z,w\in K$ we have $\val(z+w)\le ``\val(z)+\val(w)"$,
$\val(zw) = ``\val(z)\val(w)"$ and $\val^{-1}(-\infty)=0$.
In this sense $\val$ is a sub-homomorphism.
This non-Archimedean point of view on Tropical Geometry
is taken e.g. in \cite{EKL}, \cite{RGST}, \cite{Sp-diss} and \cite{SS}.
\end{remark}

\subsection{Polynomials, regular and rational functions}

There is no subtraction in the semifield $\T$, but we do not need
it to form polynomials. A tropical polynomial is a tropical sum of monomials,
for a polynomial $f$ in $n$ variables we get $f:\T^n\to\R$,
\begin{equation}\label{trop-polynomial}
f(x)=``\sum\limits_j a_jx^j"=
\operatorname{max} (a_j+<j,x>),
\end{equation}
where $x=(x_1,\dots,x_n)\in\T^n$, $j=(j_1,\dots,j_n)\in\Z^n$,
$x^j=x_1^{j_1}\dots x_n^{j_n}$, $<j,x>=j_1x_1\dots j_nx_n$
and $a_j\in\T$.

\begin{remark}
Not all monomials in $f$ are essential. It may happen that
for some element $j'$ we have
$a_{j'}+<j',x>\ \le a_j+<j,x>$
for any $j\neq j'$.
Then for any $b_{j'}\le a_j$
we have
$``\sum\limits_ja_jx^j"=``b_{j'}+\sum\limits_{j\neq j'} a_jx^j"=``\sum\limits_{j\neq j'} a_jx^j".$
Thus, the presentation of a function $f$ as a tropical polynomial is not, in general, unique.
Nevertheless, it is very close to being unique.

Note that \eqref{trop-polynomial} defines the Legendre transform of a (partially-defined)
function on $\R^n$, $-a:j\mapsto -a_j$. The function $-a$ is defined only on a finite
set, namely the set $J$ which consists of the powers of the monomials in $f$.
We can use involutivity of the Legendre transform to recover
the coefficients $a_j$ from the function $f$.
Take the Legendre transform $L_f$ of $f$ and set $\bar a_j=-L_f(j)$, $j\in\Z^n$.
It is easy to see that if $f$ were a tropical polynomial then $a_j=-\infty=0_\T$
for all but finitely many $j\in\Z^n$. Furthermore, if the function
$-a$ is convex (i.e. if it is a restriction of a convex function on $\R^n$)
then $\bar a_j=a_j$. Thus every function $f:\T^n\to\T$ obtained from a tropical
polynomial has a ``maximal" presentation as a tropical polynomial
(clearly, $\bar a_j\ge a_j$).
\end{remark}

Once we have polynomials we may form rational functions as the tropical
quotients, i.e. the differences of two polynomials.
The tropical quotient of two monomials is called a {\em $\Z$-affine} function.
Clearly a $\Z$-affine function is an affine-linear function $\R^n\to\R$
with an integer slope.

A rational function $h=``\frac fg"$ is defined for every $x\in\T^n$ such that $g(x)\neq-\infty$.
At such $x$ $h$ is called {\em finite} (more precisely, $h$ is called finite
at $x\in\T^n$ if it {\em can be presented} as $``\frac{f}{g}"$
with $g(x)\neq-\infty$).
A rational function $h$ is called {\em regular} at $x\in\T^n$
if it is finite at $x$ and
there exist a $\Z$-affine function $\phi$ finite at $x$,
an open neighborhood $U\ni x$
and a tropical polynomial $p:\T^n\to\T$ such that
$h(y)=p(y)+\phi(y)$
for all $y\in U$. (Here we consider
the Euclidean topology on $\T^n=[-\infty,+\infty)^n$.)

All functions $U\to\T$ that are restrictions of rational functions
on $\T^n$ that are regular at every point of $U$ form a semiring $\OO(U)$.
Constant functions give a natural embedding $\T\subset\OO(U)$ and
thus $\OO(U)$ is a {\em tropical algebra}. In this way we get
a sheaf $\OO$ of tropical algebras on $\T^n$ called {\em the structure sheaf}.

\section{Geometry: tropical varieties}
%piecewise-linear objects as {\em tropical varieties}}
%{\em en lieu} of algebraic varieties}

\subsection{Hypersurfaces}

To every tropical polynomial $f$ one may associate its hypersurface $V_f\subset\T^n$
that, by definition, consists of all points $x$ where $``\frac{1_{\T}}f"$ is not regular.
This is a piecewise-linear object in $\T^n$.

\begin{example}
Figures \ref{line-conic} and \ref{2cubics}
depict some hypersurfaces in $\T^2$ (i.e. planar curves).

\begin{figure}[hh]
   \includegraphics{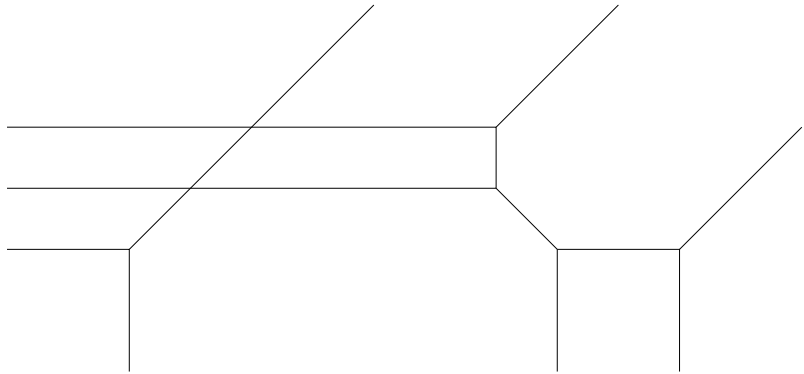}
   \caption{\label{line-conic} Tropical line and conic.}
\end{figure}

The left-hand side of Figure \ref{line-conic} is given
by the polynomial
$``1+0x+0y".$
It is a line in the tropical plane $\T^2$.
The right-hand side of Figure \ref{line-conic} is given
by the polynomial
$``10+5.5x+0x^2+8.5y+6.5y^2+4.5xy",$
it is a conic. The line and the conic here intersect at
two distinct points.

\begin{figure}[hh]
   \includegraphics{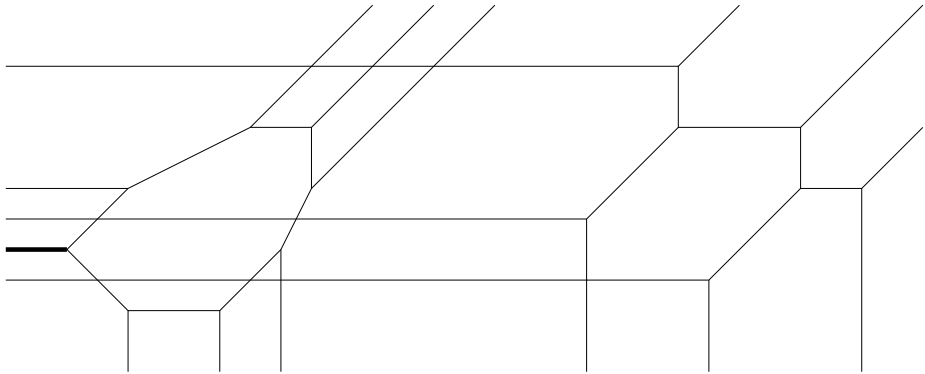}
   \caption{\label{2cubics} Two tropical cubic curves.}
\end{figure}

The left-hand side of Figure \ref{2cubics} is given
by the polynomial
\begin{equation}\label{left-cubic}
``5+4x+2.25x^2+0x^3+
4y+2.5xy+1x^2y +
3y^2+1.5xy^2+1.5y^3".
\end{equation}
The right-hand side is given
by the polynomial
$$``17.5+ 12.25x + 7x^2 + 0x^3 +
16.75y + 12xy + 5.5x^2y +
15.5y^2 + 10xy^2 + 13y^3."$$
\end{example}

It is easy to see that a tropical
hypersurface $V_f$ is a union of convex $(n-1)$-dimensional
polyhedra called the facets of $V_f$.
Each facet has {\em integer slope}, i.e. is parallel to a
hyperplane in $\R^n$ defined over $\Z$.
For every facet $P$ there
exist two monomials $a_{j_1}x^{j_1}$ and $a_{j_2}x^{j_2}$ of $f$
that are equal along $P$ and such that they are greater than
any other monomial of $f$ in the (relative) interior of $P$.
The greatest common divisor of the components of
the vector $j_2-j_1$ is called {\em the weight} of $P$.
E.g. the horisontal edge on Figure \ref{2cubics}
adjacent to the leftmost vertex has weight 2 while all
other edges has weight 1.

Intersection of any pair of facets is the face of both facets.
Furthermore, at every $(n-2)$-dimensional face $Q$ we have
the following {\em balancing property}.
\begin{property}\label{bal-hyp}
Let $P_1,\dots,P_l$
be the facets adjacent to $Q$ and let $\nu_j$
be the integer covectors whose kernels are parallel to $P_j$
and whose GCD equals to the weight of $P_j$. In addition we
require that the orientation of $\nu_1,\dots,\nu_l$ is consistent
with a choice of direction around $Q$. Then we have
$$\sum\limits_{j=1}^l\nu_j=0.$$
\end{property}

\subsection{Integer polyhedral complexes}

This balancing property may be generalized also to some
piecewise-linear polyhedral complexes $X$ of arbitrary codimension in $\T^n$.
%For simplicity we formulate this property in the absence of weights.
A {\em integer convex polyhedron} in $\T^n$ is the set defined
by a finite number of inequalities of the type
$$<j,x> \ \ \ \le \ \ \ c,$$
where $x\in\T^n$, $j\in\Z^n$ and $c\in\R$.
Here the expression $<j,x>$ stands for the scalar product.
It may happen that $<j,x>$ is an indeterminacy (for some $x\in\T\setminus\T^\times$),
then we include such $x$ into the polygon.
Equivalently, an integer convex polyhedron in $\T^n$ is
the closure of a convex polyhedron (bounded or unbounded) in $\R^n$
such that the slopes of all its faces (including the polyhedron itself)
are integers. The dimension of an integer convex polyhedron is its topological
dimension.

An integer piecewise-linear polyhedral complex $X$ of dimension $k$ in $\T^n$
is the union of a finite collection of integer convex polyhedra of dimension
$k$ called {\em the facets of $X$} such that the intersection
$\bigcap\limits_{j=1}^l P_{j}$ of any finite
number of facets is the common face of $P_j$.
We may equip the facets of $X$ with natural numbers
called the weights.

The complex $X$ is called balanced if the following property holds.
\begin{property}\label{bal-gen}
Let $Q$ be a face of dimension $k-1$ and $P_1,\dots,P_l$ be the facets
adjacent to $Q$. The affine-linear space containing $Q$ defines
a linear projection $\lambda:\R^n\to\R^{n-k+1}.$
The image $q=\lambda(Q)$ is a point while the images of $p_j=\lambda(P_j)$
are intervals in $\R^{n-k+1}$ adjacent to $q$.
Let $v_j\in\Z^n$ be the primitive integer vector parallel
to $p_j$ in the direction outgoing from $q$ multiplied by the
weight of $P_j$.
We have $$\sum\limits_{j=1}^l v_j=0.$$
\end{property}

It is easy to see that if $k=n-1$ then Properties \ref{bal-hyp} and \ref{bal-gen}
are equivalent.

%\subsection{A counterpart of the Zariski topology: open morphisms to replace open sets}
%\footnote{The idea
%to use tropical morphisms to build a reasonable topology on $\T^n$
%is due to O.Viro.}

\subsection{Contractions}

%In this paper we use Euclidean topology on $\T^n$ for all our needs.
%It is a natural topology in the tropical context
%(it can even be interpreted as a spectrum topology coming
%from tropical algebra).
%Yet sometimes
%it is convenient to consider another ``almost topology" in $\T^n$
%where the role of open sets is taken by morphisms $U\to\T^n$
%%(or, more generally, sequences of morphisms $U_k\to\dots U_1\to\T^n$)
%%rather than subsets $U\subset\T^n$.
%This correspond to
%the {Zariski topology} in the tropical set-up
%and is defined as follows.

Let $f:\T^n\to\T$ be a polynomial.
We define its {\em full graph} $$\Gamma_f\subset\T^n\times\T$$
to be the hypersurface defined by $``y+f(x)"$.
Note that $\Gamma_f$ can be obtained from
the set-theoretical graph of $f$ by attaching
the intervals $[(x,-\infty),(x,f(x))]$ for all
$x$ from the hypersurface $V_f$ (i.e. those $x$
where $``\frac 0f"$ is not regular), cf. Figure \ref{open}
for the full graphs of $``x+0"$ and $``x^2+x+1"$.

Thus, unlike the classical situation, the full graph of a map is different
from the domain of the map.
We define the
{\em principal contraction}
\begin{equation}\label{pr-open}
\delta_f:\Gamma_f\to\T^n,
\end{equation}
associated to $f$ to be the projection onto $\T^n$.
%where $D_f=\Gamma_f\cap (\T^n\times\T^\times)$
%$\Gamma_f\subset\T^n\times\T$ is the graph of $f$
%and $\delta_f$ is the projection onto $\T^n$.
%The principal open morphism $\delta_f$ extends to the {\em principal
%contraction} $$\gamma_f:\Gamma_f\to\T^n.$$

%We have $\Gamma_f,D_f\subset\T^n\times\T^\times\subset\T^{n+1}$
%and $\Gamma_f=\bar{D_f}$.
To get a general contraction one iterates this procedure.
Suppose that a contraction $\gamma:V\to\T^n$ is already defined.
%n open morphism
%which is already constructed (e.g. $\delta=\delta_f$ for some
%polynomial) and $U\subset\T^N$.
We have $V\subset\T^N$.

If $g:V\to\T$ is a regular function (in $N$ variables)
then
%we can iterate this construction by considering
we can define
%$$\Gamma_{g}=\Gamma_g\cap (U\times\T^\times)$$
%and the composition of $\delta$ and the projection $D_{g}$ to $U$
%for $\delta_{g}:D_{f,g}\to\T^n$.
%Again, we define
the full graph $\Gamma_g\subset V\times\T$
as the union of the set-theoretical graph of $g$ with
all intervals $[(x,-\infty),(x,g(x))]$ such that
$``\frac 0g"$ is not regular at $x$ (i.e. at a neighborhood of $x$ in $V$).
The map $\delta_g:\Gamma_g\to V$ is the projection onto $V$.

The map $\delta_g$ is called a principal contraction to $V$.
A general contraction is a composition of principal contractions.

Furthermore, one may associate {\em the weights} to the new facets
of $\Gamma_g$
%(and thus to the facets of $D_{g}$)
by setting them equal to the {\em orders} of the pole of $``\frac 0g"$.
%The order of the pole
%of a rational fuinction on $D_f$
%may be defined inductively,
Here we say that that the order of a pole of
a rational function is at least $n$ at $x$ if it can be locally presented
as a tropical product of $n$ rational functions that are not regular
at $x$. Then we define the order of the pole to be the largest $n$
with this property.
\begin{definition}
{\em The order of zero} of $g$ at $x$ is
the order of the pole of $``\frac 0g"$ at $x$.
\end{definition}
This definition is consistent with the definition of tropical hypersurfaces.

The facets of $\Gamma_{g}$ contained in the set-theoretical graph
of $g$ inherit their weights from the weights of the corresponding facets of $V$.

%Finite iterations of this process give all Zariski
%open morphisms $U\to\T^n$.

\begin{remark}\label{Zopen}
Contractions may be used to define counterparts of Zariski open
sets in Tropical Geometry. We define the
complements of tautological embeddings $\T^k\to\T^n$, $k\le n$,
to be Zariski open.

To define the principal open set corresponding
to a polynomial $f$ we consider the contraction
$\delta_f:\Gamma_f\to\T^n$ and take
$$D_f=\Gamma_f\cap(\T^n\times\T^\times).$$
This together with the map $\delta_f|_{D_f}$
is the principal open set associated to $f$.
\end{remark}

\begin{example}
%The embedding
%$$\T^\times\subset\T$$
%is the principal morphism for the polynomial $f(x)=x$.
%It can be considered as the result of puncturing $\T$
%at $-\infty$. Since it is an embedding, it can be considered
%as a subset of $\T$.
%
The principal open set associated to $``x+a"$
$D_{``x+a"}\to\T$
of $``x+a"$, $a\in\T^\times$,
can be interpreted as the result puncturing of $\T$ at a finite point $a\subset\T$.
%This is no longer embedding but a contraction of a ray of $D_{``x+a"}$.
Figure \ref{open} depicts $D_{``x+0"}$ and $D_{``x^2+x+1"}$.
The corresponding maps are projections onto the $x$-axis.

\begin{figure}[hh]
   \includegraphics{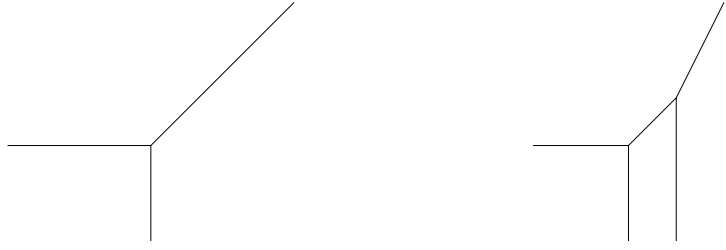}
   \caption{\label{open} Once and twice punctured affine lines $\T$.}
\end{figure}

\end{example}

\begin{proposition}\label{open-bal}
If $\delta: \T^N\supset V\to\T^n$ is a contraction then
$V\subset\T^N$ satisfies to Property \ref{bal-gen}.
%the closure of $U$ in $\T^n$ is a balanced integer piecewise-linear
%polyhedral complex.
\end{proposition}

\begin{remark}
In fact Proposition \ref{open-bal} can be used to define the
weights of the facets of $V$ and thus the orders of the zeroes
of polynomials on $V$.
\end{remark}

\subsection{Tropical varieties and tropical morphisms}

A map $\R^{N}\to\R^M$ is called {\em integer affine-linear}
if it is a composition of a linear map defined over $\Z$ and
a translation by an arbitrary vector in $\R^M$.
Such a map can be extended to a partially defined map
$\T^N\to\T^M$ by taking the closure. Note that integer affine-linear
maps leave the class of integer piecewise-linear polyhedral complexes
invariant. Furthermore, such maps take facets to facets
(at least for some presentaion of the image as
an integer piecewise-linear polyhedral complex)
and respect Property \ref{bal-gen}.

Let $X$ be a topological space together with an atlas $\{U_\alpha\}$,
$\phi_\alpha:U_\alpha\to \T^n$, $U_\alpha\subset X$.

\begin{definition}
We say that $X$ is a {\em tropical variety} of dimension $n$
if the following conditions hold.
\begin{enumerate}
\item Each $\phi_\alpha$ is a contraction to an open subset of $\T^n$.
More precisely, there is a contraction $\delta_\alpha:V_\alpha\to\T^n$,
$V_\alpha\subset\T^{N_\alpha}$ and an open embedding
$\Phi_\alpha:U_\alpha\to V_\alpha$ such that
$\phi_\alpha=\delta_\alpha\circ\Phi_\alpha$.
\item The overlapping maps $\Phi_\beta\circ\Phi_\alpha^{-1}$
are induced by the integer affine-linear maps $\R^{N_\alpha}\to\R^{N_\beta}$.
\item For every point $x\in X$ in the interior of a facet of $X$
there exists
a chart $\phi_\alpha$ such that $x\in U_\alpha$ and
$\phi_\alpha$ embeds some neighborhood of $x$ into $\T^n$.
\item There exist a finite covering of $X$ by open sets $W_j\subset X$
such that for every $j$ there exists $\alpha$ such that $W_j\subset U_\alpha$
and the closure of $\Phi_\alpha(W_j)$ in $\T^{N_\alpha}$
is contained in $\Phi_\alpha(U_\alpha)$.
%subcollection of $\{U_\alpha\}$.
%\item If $y\in \overline{\Phi_{\alpha}V_\alpha$
\end{enumerate}
\end{definition}
The last condition ensures {\em completeness} of the tropical structure on $X$.

%Let $U\subset X$ be an open set in a tropical variety.
A function on a subset $W$ of $X$ is called regular at a point $x\in W$
if it is locally a pull-back of a regular function in a neighborhood
of $\Phi_\alpha(x)\in \T^{N_\alpha}$. In this way
we get the structure sheaf $\OO_X$ of regular functions on $X$.
%\begin{proposition}
%We have a well-defined (independent on the choice of a charts)
%\end{proposition}

A point $x\in X$ is called {\em finite} if it is mapped to
$\R^{N_\alpha}\subset\T^{N_\alpha}$.
%by $\Phi_\alpha$ for $U_\alpha\ni x$.
Note that finiteness does not depend on the choice of a chart.
At finite points of $x$ we have the notion of an {\em integer tangent}
vector to $X$. This comes from tangent vectors to $\Phi_\alpha(U_\alpha)$
at $\Phi_\alpha(x)\in\R^{N_\alpha}$ after identification with the corresponding
counterparts for different charts $U_\beta\ni z$ under the differential
of the overlapping maps.

\begin{example}\label{tpn}
The space $\T^n$ is a $n$-dimensional tropical variety tautologically.
The {\em projective $n$-space $\tp^n$} is defined as the quotient
of $\T^{n+1}\setminus (-\infty,\dots,-\infty)$ by the equivalence relation
$(x_0,\dots,x_n)\sim (``\lambda x_0",\dots,``\lambda x_n")=
(\lambda+x_0,\dots,\lambda+x_n),$
where $\lambda\in\T^\times$.
%As usual, we denote the equivalence
%class of $(x_0,\dots,x_n)$ with $(x_0:\dots:x_n)$.

The space $\tp^n$ is a tropical variety since it admits (as in the classical
case) $n+1$ affine charts to $\T^n$ by dividing all coordinates by $x_j$
as long as $x_j\neq-\infty$. This is an example of a compact tropical variety.
There is a well-defined notion of a hypersurface of degree $d$ in $\tp^n$.
It is given by a homogeneous polynomial of degree $d$ in $n+1$ variables.
This polynomial can be translated to an ordinary polynomial in every affine chart
of $\tp^n$.

In  a similar way, one may construct tropicalizations of other toric varieties.
The finite part of all tropical toric varieties is $(\T^\times)^n\approx\R^n$.
\end{example}

\begin{proposition}
If $X$ and $Y$ are tropical varieties then $X\times Y$ has a natural
structure of tropical variety of dimension $\dim X+\dim Y$.
\end{proposition}

\begin{definition}\label{lin-tr-mor}
Let $f:X\to Y$ be a map between tropical varieties
(not necessarily of the same dimension).
We say that $f$ is a {\em linear tropical morphism} if for
every $x\in X$ there exist charts $U_\alpha^X\ni x$
and $U_\beta^Y\ni f(x)$ such that
$\Phi^Y_\beta\circ f\circ (\Phi_\alpha^X)^{-1}$
is induced by an integer affine-linear map
$\R^{N_\alpha}\to\R^{N_\beta}$.

The map $f$ is called a {\em regular tropical morphism}
if $\Phi^Y_\beta \circ f\circ (\Phi_\alpha^X)^{-1}$ is given
by $N_\beta$ rational functions on $\T^{N_\alpha}$ that are
regular on $\Phi_\alpha^X(U_\alpha^X)$.
\end{definition}

Clearly any linear morphism is a regular morphism.
A regular morphism $f:X\to Y$ defines a map
$\OO_Y(U)\to\OO_X(f^{-1}(U))$. This map can be interpreted
as a homomorphism of tropical algebras (defined over
the semifield $\T$).

\begin{proposition}
If $f:X\to Y$ is a linear tropical morphism then its (set-theoretical)
graph is a $(\dim X)$-dimensional tropical variety.
\end{proposition}

\subsection{Equivalence of tropical varieties}

Different tropical varieties may serve as different models
for essentially the same variety. To identify such tropical
varieties we globalize the notion of contraction that was
so far defined only for $V\subset\T^n$.

Let $f:X\to Y$ be a tropical morphism between tropical
varieties of the same dimension.
\begin{definition}
The map $f$ is called a {\em contraction} if %it is surjective
for every $y\in Y$
there exist a chart $U_\beta^Y\ni y$
with $\Phi_\beta^Y(U_\beta^Y)\subset V_\beta\subset\T^{N_\beta}$,
a contraction $\delta:W\to V_\beta$ and an isomorphism
of polyhedral complexes
$$h:f^{-1}(U_\beta^Y)\approx \delta^{-1}(\Phi_\beta^Y(U_\beta))$$
such that $\Phi_\beta^Y\circ f = \delta\circ h$.
\end{definition}

Note that a composition of contractions is again a contraction.
Contraction generates an equivalence relation on the class of
tropical varieties: tropical manifolds $X$ and $Y$ are called
equivalent if they can be connected by a sequence of contractions
or the operations inverse to contractions.

\begin{example}
The cubic curve given by $\eqref{left-cubic}$ and pictured on
the left-hand side of Figure \ref{2cubics} is equivalent to
the circle $S^1$ equipped with the tropical structure coming
from $\R/4.5\Z$ (as $\R=\T^\times$ is a tropical variety and
the translation by 4.5 is a tropical automorphism there is a
well-defined tropical structure on the quotient).

In the same time the real number 4.5 is an inner invariant of
this cubic curve. It is a tropical counterpart of the $J$-invariant
of elliptic curves.
\end{example}

It is convenient to identify equivalent tropical varieties.
This allows to present an arbitrary regular morphism $f:X\to Y$
in Definition \ref{lin-tr-mor} by a linear tropical morphism.

\begin{proposition}
If $f:X\to Y$ is a regular tropical morphism then there exists
a contraction $\delta:\tilde{X}\to X$ such that
$f\circ\delta:\tilde{X}\to Y$
is a linear tropical morphism.
\end{proposition}

\begin{example}
The map $\T\to\T$ defined by $x\mapsto f(x)=``x^2+x+1"$ is a regular morphism
which is not linear.
However, the map $f\circ\delta_f:\Gamma_f\to\T$ is a linear morphism
as it is given by the projection of the full graph $\Gamma_f$ onto the
vertical axis (cf. the left-hand side of Figure \ref{open}).
\end{example}

There is a well-defined notion of a $k$-form on a tropical variety
that is preserved by the tropical equivalence.

\begin{definition}
A $k$-form on a tropical variety $X$ is an exterior real-valued $k$-form
of the integer tangent vectors at every finite point $x\in X$
such that for every chart $U_\alpha$ this form is induced from
a (constant) linear $k$-form on $\R^{N_\alpha}$.

A $k$-form $\omega$ on $X$ is called globally defined if the following
condition holds for every non-finite point $x\in X$.
If $x\in U_\alpha$ and $\Phi_\alpha(x)\in\T^{N_\alpha}$
has its $j$th coordinate $-\infty$ then the form $\omega$ in $\R^{N_\alpha}$
vanishes on the kernel of projection onto the $j$th coordinate hyperplane.
\end{definition}

E.g. the only globally defined $k$-form on $\tp^n$ is the zero form.
But there might be globally defined $k$-forms on other compact varieties.
An easy example is provided by taking $X$ to be a {\em tropical torus},
the quotient of $\R^n$ by translation from some integer lattice $\Lambda\subset\R^n$
of rank $n$. Such $X$ is a compact tropical variety
while the globally defined $k$-forms on $X$ are in 1-1 correspondence with
constant linear $k$-forms on $\R^n$.

If $\delta:\tilde{X}\to X$ is a contraction then a globally defined
$k$-form must vanish at all vectors in the kernel of $d\delta$, thus
there is a 1-1 correspondence between forms on $\tilde{X}$ and $X$.

\section{Tropical intersection theory}

%\subsection{Subvarieties and divisors}
%A linear tropical morphism $f:Y\to X$ is called a subvariety $Y\subset X$
%if it is an embedding.

\subsection{Cycles in $X$}

The notion of an integer piecewise-linear polyhedral complex
may be extended to include not only complexes in $\T^n$, but also
complexes in an arbitrary tropical variety $X$.
We say that $B\subset X$ is an integer piecewise-linear polyhedral complex
if for every chart $U_\alpha\subset X$
there exists a integer piecewise-linear polyhedral complex
$B_\alpha\subset V_\alpha\subset\T^{N_\alpha}$
such that $\Phi_\alpha(B)\subset B_\alpha$.

As the overlapping maps preserve convex polyhedra
we have a well-defined notion of a maximal facet on $B$
and thus can consider {\em weighted} integer piecewise-linear
polyhedral complexes in $X$.

\begin{definition}
A {\em $k$-cycle} $B$ in a tropical variety $X$ is a $k$-dimensional
integer piecewise-linear polyhedral complex weighted by integer
(possibly negative) numbers
that satisfies to Property \ref{bal-gen} in every chart
$\Phi_\alpha:U_\alpha\to V_\alpha\subset\T^{N_\alpha}$ of $X$.
Accordingly, the codimension of $B$ is $n-k$.
\end{definition}

All $k$-cycles in $X$ form a group by taking unions.

\begin{remark}
In this definition we excluded $k$-cycles with boundary components.
Indeed, by our definition, every $k$-dimensional
convex polyhedron in $\T^n$ is the closure of a convex polyhedron in $\R^n$.
Thus we do not have components lying totally in $\T^n\setminus\R^n$.
While it is useful also to consider
cycles with boundary components, their intersection theory is more
elaborate.
\end{remark}

\begin{proposition}\label{push}
The image $f_*(B)\subset Y$ of a $k$-cycle $B$ under a linear tropical morphism
$f:X\to Y$ is a $k$-cycle.
\end{proposition}

In particular, a tropical subvariety is a cycle. An important example
is the {\em fundamental cycle} of the $n$-dimensional tropical variety $X$.
It is the $n$-cycle where each facet of $X$ is taken with its own weight.

\subsection{Cycle intersections}

One very useful feature of tropical varieties is the possibility to
intersect cycles there.

Let $B_1, B_2$ be two cycles of codimension $k_1$ and $k_2$
in the same tropical variety $X$.
The goal of this section is to define their {\em product-intersection} $B_1.B_2$
as a cycle of codimension $k_1+k_2$ in $X$.

%\subsection{Intersection in $\T^n$}

We start from an easier case when the ambient variety $X$ is $\T^n$.
The set-theoretical intersection $B_1\cap B_2$ is naturally
stratified by convex polyhedra that are intersections of the
(convex) faces of $B_1$ and $B_2$.
%Each open stratum of $B_1\cap B_2$
%is an intersection of two open convex polyhedra from $B_1$ and $B_2$.
It might happen that the dimension of some of these polyhedra is
greater than $n-(k_1+k_2)$.

\begin{definition}
We define $B_1.B_2$ as the closure of the union of the
strata of dimension $n-(k_1+k_2)$
in $B_1\cap B_2$ equipped with certain weights which we define as follows.
\begin{itemize}
\item Suppose that an $(n-(k_1+k_2))$-dimensional
stratum $S\subset B_1\cap B_2$ is the intersection of two
facets $F_1\subset B_1$ and $F_2\subset B_2$. Let $\Lambda_1,\Lambda_2\subset\Z^n$
be the subgroups consisting of all integer vectors parallel to
$F_1$ and $F_2$ respectively. Since by assumption $S$ is of codimension $k_1+k_2$
the sublattice $\Lambda_1+\Lambda_2\subset\Z^n$ is of finite index.
We set the weight of $S$ equal to the product of this index and the weights
of $F_1$ and $F_2$.
\item Suppose that the $(n-(k_1+k_2))$-stratum $S\subset B_1\cap B_2$
is the intersection of perhaps smaller-dimensional faces $G_1\subset B_1$
and $G_2\subset B_2$. We choose a small vector $\stackrel{\to}{v}\in\R^n$
in the generic (non-rational with non-rational projections) direction
and denote with $\tau_v:\T^n\to\T^n$ the translation by $v$ (which,
clearly, extends from $\R^n$ to $\T^n$).
The face $G_1$ is adjacent to a finite number of facets $F_1^{(\alpha)}$ of $B_1$
while the face $G_2$ is adjacent to a finite number of facets $F_2^{(\beta)}\subset B_2$.
As $v$ is chosen to be generic the facets $F_1^{(\alpha)}$ and
$\tau_v(F_2^{(\beta)})$ intersect transversely along a convex polyhedron
parallel to $S$. Thus the weight of
their intersection is already defined. We assign to $S$ the weight equal
to the sum of the weights of all such intersections (where the weight
equals to zero if $F_1^{(\alpha)}$ and $\tau_v(F_2^{(\beta)})$ are disjoint).
Proposition \ref{int} asserts that this total sum does not depend on
the choice of $v$.
\end{itemize}
\end{definition}

This definition is essentially the same as the definition of
{\em stable intersection} from \cite{RGST}, note also similarities
with the {\em Minkowski weights} from \cite{FS}.

%Let us return to the general case when the ambient variety
%for the cycles $B_1,B_2$ in an arbitrary tropical variety $X$.
%Again we define the product-intersection $B_1.B_2$ as the
%closure of the $(n-(k_1+k_2))$-dimensional strata equipped with
%some weights.
%
%To define the weights we use the contraction
%charts $\phi_\alpha$.
%It suffices to define $B_1.B_2$ in the closure
%of every facet of $X$. For every point $x\in X$
%and a choice of the facet $F$ adjacent to $x$
%there is a contraction chart
%$\phi_\alpha:U_\alpha\to\T^n$ that is an embedding at that face.
%The product-intersection $\phi_\alpha(B_1).\phi_\alpha(B_2)$
%defines $B_1.B_2$ in the closure of $F$ near $x$.
In the general case one can use in a certain way the contraction charts $\phi_\alpha : U_\alpha \to \T^n$
and the product-intersections $\phi_\alpha(B_1).\phi_\alpha(B_2)$ to define $B_1.B_2$ for an arbitrary
tropical variety $X$.

\begin{proposition}\label{int}
The product-intersection $B_1.B_2$ is well-defined.
It gives an $(n-(k_1+k_2))$-dimensional cycle in $X$.
The operation of taking the product-intersection
is commutative and associative.
\end{proposition}

\subsection{Pull-backs, deformations, linear equivalence}

Proposition \ref{push} allows us to take
the push-forward of a cycle under a linear tropical morphism.
Using the product-intersection we can also define a pull-back.

Let $f:X\to Y$ be a linear tropical morphism and $B\subset Y$
is a $k$-cycle. The product $X\times B$ is a $(\dim X+k)$-cycle
while the (set-theoretical) graph $\Gamma_f$ of $f$ is a $(\dim X)$-cycle
in $X\times Y$. Their product-intersection $(X\times B).\Gamma_f$ is
thus a $(k+\dim X-\dim Y)$-cycle in $X\times Y$.
We define the pull-back of $B$ by
$$f^*(B)=\pi^X_*((X\times B).\Gamma_f)\subset X,$$
where $\pi^X:X\times Y\to X$ is the projection onto $X$.

Any cycle $B\subset X\times Y$ can be considered as
a family of cycles in $X$. Indeed, every $y\in Y$
defines a cycle $$B_y=\pi^X_*(B.(X\times\{y\}))\subset X.$$
\begin{definition}\label{deform}
We call such a family {\em algebraic} and two cycles that
appear in the same family with a connected $B$
results of {\em deformation} of each other.
Two cycles are called {\em linearly equivalent} if they appear in the
same family with $Y=\tp^1$.
\end{definition}

\begin{proposition}
If $B_1, B_2\subset X$ are two cycles and $B'_1, B'_2\subset X$ are results
of their deformation then $B'_1.B'_2$ is a result of deformation of $B_1.B_2$.
\end{proposition}

The deformations are especially interesting in the case when $X$ is compact
as the following proposition shows. Note that a $0$-cycle in $X$ (as by
our assumption $X$ is covered by a finite collection of charts)
is a finite union of points.
We define the {\em degree} of a $0$-cycle
$B\subset X$ to be the sum of the weights of all points of $B$.

\begin{proposition}
The degree of a $0$-cycle in a compact tropical variety is a deformation invariant.
\end{proposition}

We may define the {\em intersection number} of a collection of cycles
of total codimension $n$ as the degree of their product-intersection.

\begin{example}
All hypersurfaces of the same degree in $\tp^n$ (see Example \ref{tpn})
can be obtained from each other by deformation. Furthermore, they are
linearly equivalent.
We have the tropical Bezout theorem:
the intersection number of $n$ hypersurfaces of degree $d_1,\dots,d_n$
equals to $\prod\limits_{j=1}^nd_j$, cf. \cite{S9}.

We can see the illustration to this theorem on Figures \ref{line-conic}
and \ref{2cubics}. In Figures \ref{line-conic} we have a line and a conic and they
intersect in two distinct points. The weight of each of these intersection
points is 1 as the primitive integer vectors parallel to the edges at
the points of intersection form a basis of $\Z^2$.
In Figure \ref{2cubics} we have two cubics that intersect at eight points.
Out of these eight points seven have weight 1 while one, the point
of intersection of a horizontal edge with an edge of slope 2,
has weight 2. Thus the total intersection number of these two cubics equals 9.
\end{example}

\subsection{Intersection with divisors}

Taking the product-intersection simplifies in the case when one of the cycles
is a (Cartier) divisor.

\begin{definition}\label{div}
A {\em divisor} $D\subset X$ is a finite formal linear combination
with integer coefficients
of $(\dim X-1)$-cycles given by
an open covering $\{U_\alpha\}$ and regular functions
$f_\alpha:\U_\alpha\to \T$ that defines $D\cap U_\alpha$ as its hypersurface.
\end{definition}

Not every cycle of codimension 1 in $X$ is a divisor.
\begin{example}\label{nC}
Let $X\subset \T^3$ be the hypersurface given by $``x+y+z+0"$
(see Figure \ref{non-Cartier}). Let $B\subset X$ be the line
$\{(t,t,0)\ |\ t\in\T\}$. It is a 1-cycle in the 2-dimensional
tropical variety $X$. Yet it cannot be presented as a hypersurface
near the point $(0,0,0)$ which is the vertex of $X$.
\begin{figure}[hh]
\begin{center} \includegraphics{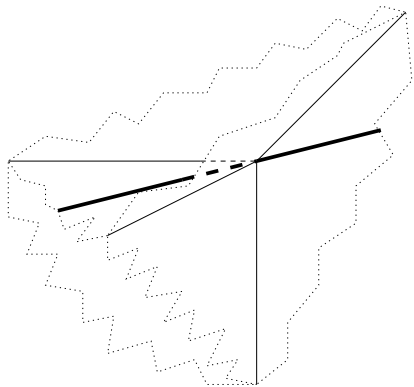} \end{center}
   \caption{\label{non-Cartier} A cycle of codimension 1 which
   is not a (Cartier) divisor.}
\end{figure}
\end{example}

If $B\subset X$ is a $k$-cycle and $D\subset X$ is a divisor
then in each $U_\alpha$ the product-intersection
$D.B$ coincides with the hypersurface defined by $f_\alpha$ on $B$
(i.e. with the points where $``\frac {1_\T}{f_\alpha}"$ is not regular
weighted by the order of poles of $``\frac {1_\T}{f_\alpha}"$).

Similarly, we may define a pull-back of a divisor $D\subset Y$ under
a regular morphism $h:X\to Y$ by taking the pull backs
of the functions $f_\alpha$. As the pull-backs $(f_\alpha)\circ h$
are regular on $h^{-1}(U_\alpha)$ they define a divisor $h^*(D)\subset X$.

%\subsection{Examples}

%\subsection{Intersections in arbitrary tropical varieties}

\section{Tropical curves}

\subsection{Tropical curves as metric graphs}

As in the classical case the easiest varieties to understand
are curves. A tropical structure on a curve (which is topologically, a graph)
can be expressed by introducing a metric. Such presentation is not unlike
the presentation of complex curves of negative Euler characteristic with
hyperbolic surfaces. For simplicity we restrict our attention to compact
tropical curves.

Recall that a {\em leaf} of a graph is an edge
adjacent to a 1-valent vertex. We call an edge which is not a leaf an
{\em finite edge}. Denote the set of all 1-valent vertices with $\ver_1$.

\begin{definition}[cf. e.g. \cite{biomath}]
A {\em metric graph} is a finite graph $\Gamma$ such that its every finite edge
is prescribed a positive real length. The length of all leaves is set to be $+\infty$.
\end{definition}

This makes $\Gamma\setminus\ver_1(\Gamma)$ a complete metric space
(equipped with an inner metric). We denote the resulting metric with $d_\Gamma$.
A homeomorphism between metric graphs is
an isomorphism if it is an isometry on $\Gamma\setminus\ver_1(\Gamma)$.
Note that a presentation of a topological space as a graph is not unique,
at our will we may introduce or erase 2-valent vertices.

\begin{proposition}
There is a 1-1 correspondence between compact tropical curves
and metric graphs.
\end{proposition}

A primitive integer tangent vector at a point of a
tropical curve (in every chart) has the unit length in the corresponding metric.

\begin{proposition}
A map $f:\Gamma\to \Gamma'$ between tropical curves is a regular
morphism if there exists a presentation of $\Gamma$ as a graph so that
for every edge $E\subset \Gamma$ there exists $n(E)\in\N\cup\{0\}$
such that $$d_{\Gamma'}(f(x),f(y))=n(E)d_\Gamma(x,y)$$
for any $x,y\in E$.
\end{proposition}

Presentation as a metric graph is a convenient way to specify
a tropical structure on a curve. The {\em genus} of a tropical
curve $\Gamma$ is its first Betti number $b_1(\Gamma)$. (This
term is justified by Proposition \ref{genus}.)

\begin{example}
Figure \ref{genus2} depicts all tropical
curves of genus 2. Varying the length of the edges
we vary the tropical structure on the curve.
In the leftmost and rightmost curves we have
three lengths of edges, $a$, $b$ and $c$, to vary. The middle curve
is an intermediate case when one of these lengths becomes zero.

\begin{figure}[hh]
   \includegraphics{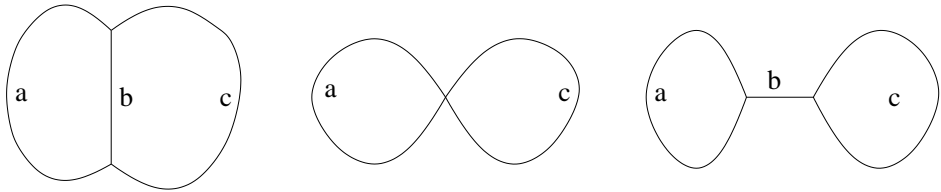}
   \caption{\label{genus2} Tropical curves of genus 2.}
\end{figure}

This shows that the space of tropical curves of genus 2 is
3-dimensional. Furthermore, one may observe that these curves
are hyperelliptic as there exists a non-trivial isometry involution
for each of the three types in the picture. (The fixed points of this isometry
are the midpoints of the edges for the leftmost picture; the midpoint of the
edges and the vertex for the center picture;
the edge connecting the vertices and the midpoints of the other edges for the
rightmost picture.)
\end{example}

\

\subsection{Jacobian varieties and the Riemann-Roch inequality}

\begin{proposition}\label{genus}
All 1-forms on a compact tropical curve $\Gamma$ form a real vector
space $\Omega(\Gamma)$ of dimension $b_1(\Gamma)$, i.e. the genus of $\Gamma$.
\end{proposition}

Note that any $0$-cycle on a curve $\Gamma$ is a divisor in the
sense of Definition \ref{div}.
%A divisor on $\Gamma$ is thus any finite combination
%f points equipped with integer weights.
The divisors of degree 0 form a group $\Div_0$ while
the divisors of degree $d$ form a homogeneous space $\Div_d$
over $\Div_0$. We set the {\em Picard group} $\pic_0$ by
taking the quotient group of $\Div_0$ by the linear equivalence
(see Definition \ref{deform}), similarly we define $\pic_d$.

A divisor is called {\em principal} if it is linearly equivalent to zero.
Clearly, the degree of a principal divisor is 0.
The classical Mittag-Leffler problem to determine whether a divisor
of degree 0 is principal is answered by the Abel-Jacobi theorem.
A similar answer exists also in the tropical set-up.

Note that given a 1-form $\omega$ and a path $\gamma:[0,1]\to\Gamma$
there is a well-defined integral $\int\limits_\gamma\omega$.
Clearly, the value of this integral depends only on the
relative homology class of $\gamma$.
Let $\Omega^*(\Gamma)$ be the (real) dual vector space to $\Omega(\Gamma)$.
Its dimension is equal to the genus $g$ of $\Gamma$ by Proposition \ref{genus}.
%Thus we may extend such integration to
Each element $a\in H_1(\Gamma;\Z)$ determines a point of $\Omega^*(\Gamma)$.
This point is given by the functional
$$\Lambda(a):\omega\mapsto\int\limits_a\omega.$$

%\begin{definition}
The {\em Jacobian} of a tropical curve $\Gamma$ is defined by
$$\Jac(\Gamma)=\Omega^*(\Gamma)/\Lambda(H_1(\Gamma;\Z)).$$
%\end{definition}
The Jacobian is an example of a tropical torus, it is homeomorphic
to $(S^1)^g$ and carries the structure of a tropical variety (which
depends on the lattice $\Lambda(H_1(\Gamma;\Z))$.

To define the {\em Abel-Jacobi} map
\begin{equation}\label{AJ}
\alpha:\Pic_0(\Gamma)\to\Jac(\Gamma)
\end{equation}
we take any 1-chain $C$ whose boundary is a divisor in the equivalence class
$p\in\Pic_0$ and set $\alpha(p)$ to be the functional
$\omega\mapsto \int\limits_C \omega$.

\begin{proposition}[Tropical Abel-Jacobi Theorem]
The Abel-Jacobi map \eqref{AJ} is a well-defined bijection.
\end{proposition}

This gives the structure of a tropical variety on $\Pic_0$ as well as
on $\Pic_d$ (which is a homogeneous space over $\Pic_0$).
We have the tautological map
$$\Sym^d(\Gamma)\to\Pic_d(\Gamma)$$
defined by taking the equivalence class.

As in the classical case this map is especially interesting in the
case $d=g-1$. The image of this map in this case (as well as all
its translates in $\Pic_0=\Jac$ is called the {\em $\Theta$-divisor}.
There is a tropical counterpart of the Riemann theorem stating
that the $\Theta$-divisor is given by a $\theta$-function, see \cite{MZ}.

\begin{figure}[h]
\begin{center}
   \includegraphics{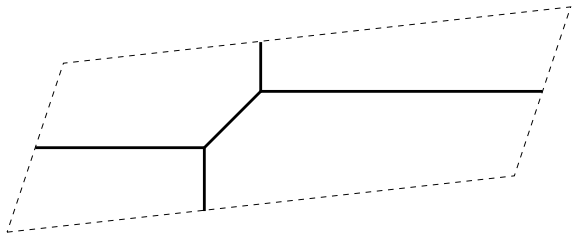}
   \end{center}
   \caption{\label{theta} $\Theta$-divisor for genus 2.}
\end{figure}

Figure \ref{theta} sketches the $\Theta$-divisor in the case of genus 2.
In this case it is isomorphic to the curve $\Gamma$ itself.
The Jacobian on Figure \ref{theta} is obtained by identifying
the opposite sides of the dashed parallelogram.

The divisor is called {\em effective} if the weights at all its points
are positive. The Riemann-Roch theorem allows to find the dimension
of all effective divisors linearly equivalent to a given one.
(Alternatively, this number plus one can be interpreted as the dimension
of the space of sections of the line bundle corresponding to
a given divisor.) Tropically, we have the Riemann-Roch theorem
in form of an inequality.

\begin{proposition}[The tropical Riemann-Roch inequality]\label{trRR}
The dimension of the space of effective divisors in the equivalence
class $p\in\Pic_d$ is at least $d-g$.
\end{proposition}

\

\ignore{
\subsection{Tropical curves in $\R^n$}

Tropical toric varieties and, in particular, $\tp^n$ give some
simple examples of tropical varieties. Their finite part is $\R^n$.
A curve in a toric variety which does not have components lying
totally in the boundary divisor can be thought of as a punctured
tropical curve in $\R^n$ (cf. Remark \ref{Zopen}).

Consider a 1-cycle $B$ in the tropical variety $\R^n=(\T^\times)^n$
such that all the weights at its edges are positive.
One can present $B$ as an image of the fundamental cycle
for a linear tropical morphism $$h:\Gamma\to\R^n$$ for
a (non-compact) tropical curve $\Gamma$. Furthermore,
we may assume that $h$ is topologically an immersion
(we may present an edge of weight $w$ as an image of
a single edge of $\Gamma$ after stretching $w$ times).
%It is possible to compactify $h$ to a map from
%The deformation type of $B$ is determined by the direction and the
%weight of its unbounded edges.

Let $x$ be the number of punctures (i.e. the number of ends) of $\Gamma$
and $g$ be its genus.
If we deform both $h$ and the tropical structure on $\Gamma$ one
may compute the dimension, or the {\em nember of degrees of freedom}
of such deformations using Proposition \ref{trRR}.
This dimension is not greater than
\begin{equation}\label{rr}
x+(n-3)(1-g).
\end{equation}
}

\subsection{The canonical class,
regular and superabundant curves in $X$}

We have the Chern classes for compact tropical varieties which
are easy to compute (as a tropical variety is already parallelized
by the integer affine structure on its facets).

Let $X$ be an $n$-dimensional compact tropical
variety. We have the natural stratification of the faces of $X$
by their dimension. Furthermore, there is the {\em boundary
stratification} of $X$. We say that a face $F\subset X$
is in the $k$th boundary stratum if there is a chart $U_\alpha\ni x$
where $\Phi_\alpha(F)$ is contained in an intersection of $k$
out of $N_\alpha$ tropical coordinate hyperplanes of $\T^{N_\alpha}$
and not contained in the intersection of any $(k+1)$
coordinate hyperplanes.

A face in the 0th boundary stratum is called a {\em finite} face.
Clearly a face in the $k$th boundary stratum
has the dimension at most $(n-k)$.

The $k$th Chern class of $X$ is an $(n-k)$-cycle that is
the linear combination of all $(n-k)$-dimensional faces of $X$.
Here the strata of boundary codimension $k$ are taken with the
weight equal to $(-1)^k$ times some positive number depending
only on the local geometry near a point $x$ in the relative interior of $F$
(this weight can be computed inductively from a local presentation
of a neighborhood of $x$ by a contraction  $\delta:V\to\T^n$).

Here we give the computation of weights only in the case of $-c_1$.
This is the $(n-1)$-cycle in $X$ that consists of all finite $(n-1)$-faces $F$
taken with the weight equal to the number of adjacent facets to $F$
minus 2 and
all $(n-1)$-faces in the 1-boundary stratum taken with the weight $-1$.

\begin{definition}
The {\em canonical class} $K$ is the $(n-1)$-cycle in $X$ equal to $-c_1$.
\end{definition}

A {\em parameterized tropical curve} in $X$
is a linear tropical morphism $h:\Gamma\to X$.
One may use Proposition \ref{trRR} to compute
the dimension in which the map $h$ varies (we allow to deform both $h$
and the tropical structure on $\Gamma$).

%\begin{example}
The adjunction formula exists also in the tropical geometry.
As in this talk we have not defined tropical vector bundles
in general and, in particular, the normal bundles, we state it only
in the case when $h$ is an embedding and $X$ is a compact surface.
We have
\begin{equation}\label{adj}
h_*([\Gamma]).h_*([\Gamma])=g-2-K.h_*([\Gamma]).
\end{equation}
%\end{example}
Here $[\Gamma]$ stands for the fundamental 1-cycle in $\Gamma$.
%The formula \ref{adj} just me
\begin{example}
Consider the compactification of Example \ref{nC}. Let $X\subset\tp^3$ be the closure
of the hypersurface $x+y+z+0$ and $h:\tp^1\to X$ be the linear tropical
morphism whose image is the closure of $B$.
The two infinite ends of $B$ contribute $+1$ each to its self-intersection, thus
the self-intersection contribution of $B$ at $(0,0,0)$ is $-1$ (see Figure
\ref{non-Cartier}).
\end{example}

\begin{proposition}
The map $h:\Gamma\to X$ varies in a family of dimension at least
$$\operatorname{\vdim_h}=K.h_*([\Gamma])+(n-3)(1-g)$$
if we allow to deform both $h$ and $\Gamma$.
The number $\vdim$ is called the {\em virtual dimension}
of the deformations of $h$.
\end{proposition}

\begin{definition}
The parameterized curve $h$ is called {\em regular} if the local dimension of
all its deformations is equal to  $\vdim_h$. Otherwise,
if it is strictly greater than $\vdim_h$, the curve $h$ is called
{\em superabundant}.
\end{definition}

E.g. the curve $h$ will necessarily be superabundant if $h(\Gamma)$ contains
a loop that is contained in a proper affine-linear subspace of a facet of $X$.

\subsection{Tropical moduli spaces}

Let $\Gamma$ be a compact tropical curve of genus $g$.
Let us choose $k$ distinct points $x_1,\dots,x_k\in \Gamma$.
We call $x_j$ the {\em marked points}.
By replacing $\Gamma$
with an equivalent tropical curve if needed,
we may assume that $\Gamma$ has exactly $k$ 1-valent
vertices which coincide with $x_1,\dots,x_k$.
Such presentation exists unless $g=0$ and $k<2$
and it is unique up to isomorphism. Thus we may restrict
our attention only to such models.

We denote with $\M_{g,k}$ the space of all tropical curves of genus $g$
with $k$ marked points up to equivalence.
\begin{figure}[hh]
\begin{center}
   \includegraphics{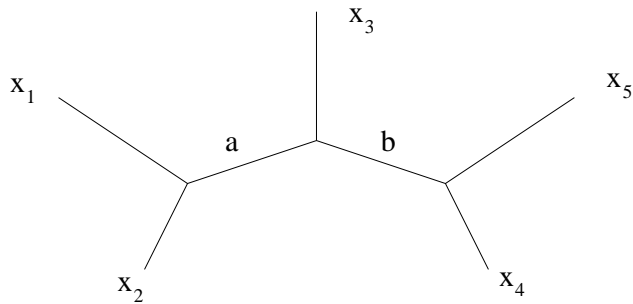}
\end{center}
   \caption{\label{mpoints} A rational curve with 5 marked points.}
\end{figure}

A tropical curve with marked points is determined by its combinatorial
type and the lengths of its finite edges. E.g. the tropical
curve from Figure \ref{mpoints} is determined by the lengths $a$ and $b$.
The length of the finite edges can be used to enhance $\M_{g,k}$ with
a tropical structure. The only problem is that we may only identify the
edges within the same combinatorial type of the curve.
%while for the tropical structure we also need

One can easily solve this problem if $g=0$.
Let us introduce a global function $Z_{x_i,x_j}$ on $\M_{0,k}$
for any pair of marked
points $\{x_i,x_j\}$. We set $Z_{x_i,x_j}$
equal to the total length of finite edges between $x_i$ and $x_j$.
E.g. for the curve $\Gamma\in\M_{0,5}$ from Figure \ref{mpoints}
we have $Z_{x_1,x_2}(\Gamma)=0$, $Z_{x_1,x_3}(\Gamma)=a$
and $Z_{x_1,x_5}(\Gamma)=a+b$.

\begin{proposition}
The functions $Z_{x_i,x_j}$ define an embedding
$$\M_{0,k}\subset\R^{\frac{k(k-1)}{2}}.$$
Its image is a tropical subvariety.
% in $\R^{\frac{k(k-1)}{2}}$.
\end{proposition}

\begin{remark}
This embedding is related to the Pl\"ucker coordinates on the
Grassmannian $G_{2,k}$, see \cite{SS} for a tropical version.
\end{remark}

Thus we see that $\M_{0,k}$ has a natural structure of a tropical variety.
We may compactify $\M_{0,k}$ by allowing the lengths of
some (or all) finite edges of $\Gamma$ to take
the value equal to $+\infty$. Such a generalized curve
splits to several components: the finite points of each such component
are within finite distance from each other.
\begin{proposition}
The resulting space $\overline\M_{0,k}$ is a compact tropical variety.
\end{proposition}
This compactification is a tropical counterpart of the Deligne-Mumford
compactification.

If $g>0$ one may use similar arguments to show that $\M_{g,k}$ is
a tropical orbifold which can be compactified to a compact
tropical orbifold $\overline\M_{g,k}$.

\subsection{Stable curves in $X$}

As in the classical case we call a parameterized
tropical curve $h:\Gamma\to X$
in a compact $n$-dimensional tropical variety $X$ with compact $\Gamma$
{\em stable} if there are no infinitesimal automorphisms of $h$,
i.e. if the number of isomorphisms $\Phi$ of $\Gamma$ such
that $h=h\circ\Phi$ is finite.

All deformations of a regular stable curve $h$ locally
form a tropical variety of dimension $K.h_*([\Gamma])+(n-3)(1-g)$.
Let us fix some class $\beta$ of stable curves to $X$ that
is closed with respect to regular curves deformations.
The class $\beta$ can be given e.g. by prescribing the intersection
number with all $(n-1)$-cycles in $X$.

Denote with $\M_{g,k}^\beta(X)$
the space of stable curves of genus $g$ with $k$ marked points
in the class.
In many cases the space $\M_{g,k}^\beta(X)$
can be compactified to a compact tropical variety $\overline\M_{g,k}^\beta(X)$.
This holds for instance if $g=0$, $X=\tp^n$ and $\beta$ is formed by curves
of degree $d$. Another instance is if $X=\tp^2$, there are no restrictions
on $g$ and $\beta$ is formed by topological immersions of degree $d$.
These are the two principal cases for our enumerative applications.
More generally we may assume that $X$ is any
compact toric variety, but then we have to impose an additional
constraint on $\beta$ that it does not contain curves with components
lying totally in the boundary divisors of $X$.

Furthermore, in these cases we have the {\em evaluation} map
$$\ev_j:\overline\M_{g,k}^\beta(X)\to X,$$ $\ev_j(h)=h(x_j)$
as well as the maps $\ft:\overline\M_{g,k}^\beta(X)\to\overline\M_{g,k},$
$\ft(h)=\Gamma$, and the maps
$\pi_j:\overline\M_{g,k}^\beta(X)\to\overline\M_{g,k-1}^\beta(X)$
``forgetting" the marked point $x_j$.
These maps are linear tropical morphisms.

This allows to set up a tropical framework for the Gromov-Witten
theory. E.g. given a collection of cycles in $X$ we may take their
pull-backs and then take their intersection number in $\overline\M_{g,k}^\beta(X)$.

Many reasonings in the Gromov-Witten theory can be literally repeated
in this tropical set-up. A good example is the WDVV-relation, cf. \cite{GM}.
As in the classical case the stable curves in
$\overline\M_{0,k}^\beta(X)\setminus\M_{0,k}^\beta(X)$
must consist of several components.

\begin{remark}
Moduli spaces of higher-dimensional tropical varieties
is a very interesting, but much more difficult subject.
Already the tropical K3-surfaces form a very
sophisticated geometric object, see \cite{KS-K3} and \cite{GW}.
%for the Differential Geometry viewpoint)
%for its description in the case of K3-surfaces.
\end{remark}

\section{Tropical curves in $R^n$, their phases and amoebas}

%\subsection{Amoebas and their degeneration}

Let $V\subset\tor$ be an algebraic variety.
Its {\em amoeba} (see \cite{GKZ}) is the set
$\am=\Log(V)\subset\R^n,$
where $\Log(z_1,\dots,z_n)=(\log|z_1|,\dots,|z_n|)$.
Similarly we may consider the map $$\Log_t:\tor\to\R^n$$
corresponding to taking the logarithm with the base $t>1$.

Amoebas themselves have proved to be a very useful tool
in several areas of Mathematics, see e.g.  \cite{EKL},
\cite{Mi-am}, \cite{Mi-pp}, \cite{MiRu}, \cite{PR}, \cite{PST}.
However, for the purposes
of this talk we only use them
as an intermediate link between the classical and tropical geometries.

\begin{definition}
The curve $h:\Gamma\to\R^n$ is called {\em classically realizable}
if there exist a small regular neighborhood $U\supset B$ in $\R^n$,
a retraction $\rho:U\to B$, a regular family of holomorphic maps
$H_t:C_t\to\tor$ for a family of Riemann surfaces $C_t$
defined for all sufficiently large positive $t>>1$ and
smooth maps $\lambda_t:C_t\to\Gamma$ such that
\begin{itemize}
\item
$h\circ\lambda_t=\rho\circ\Log_t\circ H_t.$
\item The genus of $C_t$ coincides with the genus of $\Gamma$.
\item The number of punctures of $C_t$ coincides with the number
of ends of $\Gamma$.
\end{itemize}
The family $H_t:C_t\to\tor$ is called an {\em approximating family} of $h$.
\end{definition}

%Let $E\subset\Gamma$ be a set containing a point in the interior
%of every finite edge.
\begin{proposition}
If $h:\Gamma\to\R^n$ is a tropical curve approximated
by $H_t$ then for a sufficiently large $t$ the inverse image $\lambda_t^{-1}(p)$
is a smooth circle for every $p$ inside an edge of $\Gamma$ while
the inverse image $\lambda_t^{-1}(W)$ is diffeomorphic to a sphere with $u$ punctures
for a small connected neighborhood $W$ of a vertex of valence $u$.

In particular, if $\Gamma$ is 3-valent then $\lambda_t$ defines a pair-of-pants
decomposition. In turn, this pair-of-pants decomposition determines
a point in the boundary of the classical Deligne-Mumford space $\overline\M_{g,k}^\C$.
This point is the limit of the Riemann surfaces $C_t$.
\end{proposition}
See \cite{Mi-pp} for a generalization of the pair-of-pants decomposition
for the case of higher-dimensional hypersurfaces .

\begin{remark}
Classical realizability of tropical varieties in $\R^n$ is
closely related to their presentation by {\em non-Archimedean amoebas},
the images of algebraic varieties in $(K^\times)^n$ under the coordinatewise
valuations (as defined by Kapranov \cite{Ka}).
Here $K$ is an algebraically closed field with
a non-Archimedean valuation $\val:K^\times\to\R$. See \cite{Sp-diss}
for an account of what is known on such presentations.
\end{remark}

To formulate the realizability theorem in full generality
we need to define the {\em phases} for $\Gamma$.
We start from their definition in a model case.

Let $\Gamma_k\subset\R^k$ be the tropical curve consisting
of $(k+1)$ rays emanating from $0\in\R^k$ in the directions
$(-1,\dots,0),\dots,(0,\dots,-1)$ and $(1,\dots,1)$.
The tautological embedding $\Gamma_k\subset\R^k$ is
easily realizable. For an approximating family we
can take
$H_t=L_k\subset\tor\subset\cp^k$, where $L_k$ is a
line with $(k+1)$ ends in $\tor$.
The choice of $L_k$ up to a multiplication by a point of $(\R_+)^k$
is called the phase of $\Gamma_k$.
%Let
%$\Arg:\tor\to (S^1)^n$ be the map defined by taking the argument coodinatewise.
%The phase of $\Gamma_n$ is the
%argument projection

Locally, near any point $x\in\Gamma$ the map $h$ coincides
with the map $$\Gamma_k\subset\R^k\stackrel{A+c}\to\R^n$$
near $0\in\Gamma_k$ for a linear map $A$ defined over $\Z$
and $c\in\R^n$,
where $(k+1)$ is the valence of $x$.
The linear map $A$ can be exponentiated to
a multiplicative linear map $a:(\C^\times)^k\to\tor$.
The phase $\sigma_U$ of $\Gamma$ at a small neighborhood $U\ni x$
is defined to be the equivalence class of $\xi a(L_k)\subset\tor$, $\xi\in\tor$,
up to multiplication by an element of $(\R_+)^n$.
%The phase of $\Gamma$ is a collection of compatible local
%phases in an open covering of $\Gamma$.
Two local phases are called {\em compatible} if they agree on
the intersection of the neighborhoods.

We say that $H_t$ approximates $\Gamma$ {\em with the local
phase $\sigma_U$} if $H_t(C_t)\cap\Log_t^{-1}(U)$ converges (in the Hausdorff metric on compacts in $\tor$)
to $\xi t^C a(L_k)\in\sigma_U$.

\begin{thm}\label{realize}
Any regular tropical curve $h:\Gamma\to \R^n$ equipped
with any compatible system of phases
is classically realizable.
\end{thm}
Cf. \cite{Mi-R2} and \cite{Sh},
for the special case $n=2$ with no restriction on the genus,
cf. \cite{Mi-Q} and \cite{NS} for the special case
of genus 0 with no restriction on $n$.
The proof in the general case (though in a somewhat different language)
is contained in \cite{B} (cf. also \cite{Sp-diss}).

\begin{remark}
Regularity is a necessary condition in Theorem \ref{realize}, it is easy to construct
an example of a non-realizable superabundant curve even for a topological
immersion of an elliptic curve in $\R^3$, see e.g. \cite{Mi-GV}.
\end{remark}

\section{Applications}

One of the greatest advantage of Tropical Geometry is
that most classical problems become much simpler after their tropicalization
(if such a tropicalization exists!).
This simplicity comes from the piecewise-linear nature of the tropical objects.
Indeed, once we fix the combinatorial type of the data,
a tropical problem becomes linear. E.g. if there are finitely many
solutions to a problem then in every combinatorial
type we will have a unique solution or no solutions at all.

This allows one to find (at least) an algorithmic answer to
a tropical problem. Sometimes one can show that the answer
to a classical problem and its tropicalization must coincide
and obtain the answer to the classical problem in this way.
In this last section we list some examples of such problems.

\subsection{Complex Geometry}

Let $g\ge 0$ and $d\ge 1$ be integer numbers. Fix a collection
${\mathcal Z}=\{z_j\}_{j=1}^{3d-1+g}$ of points in $\cp^2$ in general position.
There are finitely many holomorphic curves of genus $g$
and degree $d$ passing through ${\mathcal Z}$.
%The condition
%of passing through each $z_j$ defines a divisor in the space
%of all curves of genus $g$ and degree $d$ in $\cp^2$. These divisors
%intersect transversally at their smooth points if ${\mathcal Z}$ is generic.
Let $N^\C_{g,d}$ be their number.
% of holomorphic curves

Such set-up can be almost literally repeated in the tropical framework.
Fix a collection
${\mathcal X}=\{x_j\}_{j=1}^{3d-1+g}$ of points in $\tp^2$ in general position (see
\cite{Mi-R2}).
Again, it can be shown that there are finitely many tropical curves
$h:\Gamma\to\tp^n$ of genus $g$
and degree $d$ passing through $x_j$.
But in the tropical set-up these curves come with natural
positive integer multiplicities $m(h)$ not necessarily equal to 1.

Tropical curves of genus $g$ and degree $d$ that pass through $x_j\in\tp^2$
form a codimension 1 cycle in $\overline{\M}^d_{g,d}(\tp^2)$.
We set $m(h)$ to be the weight of their product-intersection
at $h$ (in other words this is their local intersection number in $h$).
%The curves passing through ${\mathcal X}$ are their intersectio
%If ${\mathcal X}$ is generic thenhese cycles intersect in a finite
%number of points and
Let $N^\T_{g,d}$ be the number of the tropical curves of degree $d$
and genus $g$ passing through ${\mathcal X}$ counted
with multiplicity $m$.
%\begin{remark}
%This number $m(h)$ is easily computable just by looking
%at the image 1-cycle $B\subset\R^2\subset\tp^2$.
%At every 3-valent vertex $v$ of $B$ we define
%the local multiplicity as $m_v(h)=|w\wedge w'|$ where $w$ and
%$w'$ are the primitive integer vectors parallel to two (out of three)
%adjacent edges at $v$
%\end{remark}

\begin{thm}[\cite{Mi-R2}]\label{cp2}
$$N^\C_{g,d}=N^\T_{g,d}.$$
\end{thm}

Of course, there are well-known ways to compute
$N^\C_{g,d}$ (see \cite{KM}, \cite{CH}) as they coincide with
the Gromov-Witten invariants of $\cp^2$.
Yet Theorem \ref{cp2} gives another simple and visual way to do it,
see \cite{Mi-R2} for details.

\begin{example}\label{deg3comp}
Figure \ref{9curves} depicts rational cubic curves in $\tp^2$
passing through a generic configuration ${\mathcal X}$ of 8 points.
There are nine curves.  Eight of them have multiplicity 1. The remaining
one, namely the rightmost in the middle row, has multiplicity 4
(it has an edge of weight 2 shown by a bold line in the picture).
Thus the total number of tropical curves is 12.
\begin{figure}[hh]
\begin{center}
   \includegraphics{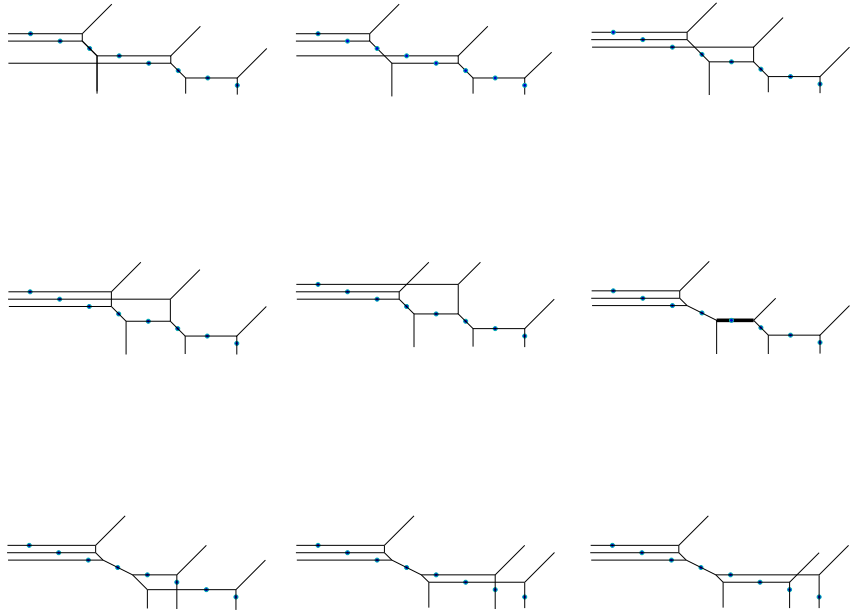}
\end{center}
   \caption{\label{9curves} Rational tropical cubics via 8 generic points in the plane.}
\end{figure}

Should we choose a different generic configuration ${\mathcal X}$
the number of tropical curves could be different. But their total
number counted with multiplicities is invariant. E.g. there
exists a different choice of ${\mathcal X}$ where we have ten tropical
curves, nine of multiplicity 1 and one of multiplicity 3. In fact,
no other partition of 12 can appear by Example \ref{rdeg3}.
\end{example}

\begin{remark}
In \cite{Mi-R2} there is also a version of Theorem \ref{cp2}
for other toric surfaces. Note, that if the ambient toric surface is not Fano
then there is a difference between counting irreducible curves
in a given homology class and computing the corresponding Gromov-Witten
number (as the latter also has a contribution from curves with some
boundary divisors as components). Tropical Geometry has an advantage
of giving a way to compute the number of irreducible curves directly
without having to deal with those extra components.

The same advantage allows to apply Tropical Geometry for giving
an algorithmic answer to another
classical problem of the Complex Geometry, computation of Zeuthen's numbers,
see \cite{M-Zeut}. These are the numbers of curves of degree $d$ and genus $g$
that pass through a collection of generic points and tangent to a collection
of generic lines in $\cp^2$. Here the total number of points and lines in
the configuration is $3d-1+g$. These Zeuthen's numbers also turn out to
be equal to the corresponding tropical numbers, while the latter can
be computed by a finite (though quite extensive for large genus) algorithm.
As far as the author knows such computation in general is not currently
accessible by non-tropical techniques.
\end{remark}
%First we note that $m(h)$ are easy to compute. For a generic ${\mathcal X}$
%the corresponding graph $\Gamma$ must be 3-valent. At any vertex $v\in\Gamma$
%we define

\subsection{Real Geometry}

The number of real curves of degree $d$ and geometric genus $g$
passing via a collection ${\mathcal Z}=\{z_j\}_{j=1}^{3d-1+g}$ of points in $\rp^2$
depends on ${\mathcal Z}$ even if we choose it to be generic. This is the
feature of $\R$ as a non-algebraically closed field.
Yet, as it was suggested in \cite{We},
one may prescribe a {\em sign} $\pm 1$ to every such curve so that
the sum $W_d$ of these signs becomes invariant in the special case of $g=0$.
The number $W_d$ is called the {\em Welschinger number}.

If ${\mathcal Z}$ is a generic configuration of $3d-1$ points in $\rp^2$
then any real rational curve $\R C$ of degree $d$ passing through ${\mathcal Z}$
is a nodal curve, i.e. all singularities of $\R C$ are non-degenerate double points.
Over $\R$ there are two types of such nodes: the {\em hyperbolic} node, corresponding
to the intersection of two real branches (given in local coordinates by $x^2-y^2=0$)
and the {\em elliptic} node, corresponding to
the intersection of two complex-conjugate branches
(given in local coordinates by $x^2+y^2=0$).
The sign of $\R C$ is defined as $(-1)^{e(C)}$, where $e(C)$
is the number of elliptic nodes of $C$.

The Welschinger number also has a tropical counterpart.
To a tropical curve $h:\Gamma\to\tp^2$ of degree $d$ and genus $g$
passing through a configuration ${\mathcal X}$ of $3d-1+g$ points
we associates its {\em real multiplicity} $m^\R(h)$ which is $\pm 1$ or $0$.
If $h(\Gamma)$ has an edge of even multiplicity then $m^\R(h)=0$.
Otherwise we define the local real multiplicity $m_v^\R(h)$ for
a vertex $v\in\Gamma$ to be $(-1)^{e_v}$, where $e_v$
is the number of integer points in the interior of the lattice triangle
such that its sides are perpendicular to the edges adjacent to $v$
and of integer length equal to the weight of that edge.
Then we define $m^\R(h)=\prod\limits_v m^\R_v(h)$.
Let $W^\T_{g,d}$ be the corresponding tropical number.

\begin{thm}[\cite{Mi-R2}]\label{wel}
$$W_d=W^\T_{g,d}.$$
\end{thm}

This theorem is the only currently known way to compute
the Welschinger numbers, see \cite{IKS}.

\begin{example}\label{rdeg3}
Let us revisit Example \ref{deg3comp}. The real multiplicities
of eight out of the nine curve of Figure \ref{9curves} is 1 while
the real multiplicity of the remaining one is 0. Thus we have $W_3=8$.
For another choice of a generic configuration of 8 points in $\tp^2$
mentioned in Example \ref{deg3comp} we get 9 curves with $m^\R(h)=+1$
and one with $m^\R(h)=-1$.

Note that we always have $m^\R(h)=+1$ if $m(h)=1$;
$m^\R(h)=+1$ if $m(h)=1$ if the multiplicity is 1;
and $m^\R(h)=0$ if $m(h)$ is even. Note also that we may never get $m(h)=2$.
Since the sum of the multiplicities has to be 12 and the sum of real
multiplicities has to be 8 the partitions $12=8+4=9+3$ are the only
two possible partitions. In particular, there does not exist a configuration
of 8 points in $\tp^2$ such that the 12 tropical cubics
%passing through them
will all be distinct.
\end{example}

\begin{remark}
There is a 3-dimensional version of the number $W_d$
which is the number of real rational curves in $\rp^3$
of degree $d$ passing through a generic configuration of $2d$ points
taken with certain signs, see \cite{We3}. These numbers
also have tropical counterparts and Theorem \ref{wel}
extends to the 3-dimensional case providing a way to compute
the real geometry numbers, see \cite{Mi-Q}.
\end{remark}

As the last application of tropical geometry we would like to mention
extending the patchworking \cite{Vi} to curves in real toric
varieties of higher dimensions by using Theorem \ref{realize}
with real phases.
%(that are determined in a 2-to-1 fashion
%by a point in $(\Z_2)^n$ near a point on the edge of $\Gamma$).

\begin{figure}[hh]
\begin{center}
   \includegraphics{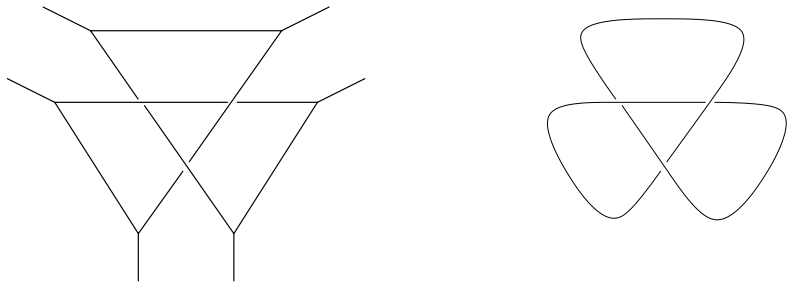}
\end{center}
   \caption{\label{knot} Real algebraic knots from tropical curves in $\R^3$.}
\end{figure}

\begin{example}
Let $K\subset\R^3$ be a knot presented as an embedded
piecewise-linear circle made with $k$ straight intervals.
Then there exists an algebraic curve in $(\R^\times)^3$
with a unique closed component in $(\R_+)^3$
isotopic to $K$ such that its complexification
has genus 1 and is punctured $k$ times. E.g. a trefoil may be presented
by an elliptic curve with 6 punctures in $(\C^\times)^3$, see Figure \ref{knot}.

To deduce this statement we perturb the broken line $K$ to make
the slopes of its intervals rational. Then we add
an extra rays at every corner to get
a tropical curve $\Gamma\supset K$. Finally we choose real phases for $\Gamma$ so that the
phases of all its bounded edges include the positive quadrant.
\end{example}

\ignore{

\section{Deadline and where to send the file}

We ask you to send your contribution (tex, dvi, pdf as well as eps files for illustrations)
\textbf{no later} than December 31, 2005 by ftp to the following account:

\smallskip\noindent
Hostname: ftp.ems-ph.org

\smallskip\noindent
Username: icm2006

\smallskip\noindent
Passwort: pfpw4eph

\smallskip\noindent
Folder: submissions

\bigskip\noindent
This would help us for timely production.

}


\frenchspacing
\begin{thebibliography}{99}

\bibitem{biomath} Billera, L. J., Holmes, S. P., Vogtmann, K.,
Geometry of the space of phylogenetic trees.
Advances in Applied Math 27, 733-767 (2001).

\bibitem{B} Bourgeois, F., A Morse-Bott approach to contact homology.
Dissertation, Stanford University 2002.

\bibitem{CH} Caporaso, L., Harris, J.,
{Counting plane curves of any genus}.
{\em Invent. Math.} {\bf 131} (1998), 345-392.

\bibitem{EKL} Einsiedler, M., Kapranov, M., Lind, D.,
Non-archimedean amoebas and tropical varieties,
http://arxiv.org/abs/math.AG/0408311.

%\bibitem{EGH}

\bibitem{Fu} Fukaya, K.,
Multivalued Morse theory, asymptotic analysis and mirror symmetry.
In {\em Graphs and patterns in mathematics and theoretical physics,}
Proc. Sympos. Pure Math., 73, Amer. Math. Soc., Providence, RI, 2005,
 205--278.

\bibitem{FS} Fulton, W., Sturmfels, B.,
Intersection theory on toric varieties.
{\em Topology}  {\bf 36}  (1997), 335--353.

\bibitem{GM} Gathmann, A., Markwig, H.,
Kontsevich's formula and the WDVV equations in tropical geometry,
http://arxiv.org/abs/math.AG/0509628.

\bibitem{GKZ} Gelfand, I. M., Kapranov, M. M., Zelevinsky, A. V.,
{\em Discriminants, resultants, and multidimensional determinants.}
Birkh\"auser Boston, Inc., Boston, MA, 1994.

\bibitem{GW} Gross, M., Wilson, P. M. H.,
Large complex structure limits of $K3$ surfaces.
{\em J. Differential Geom.}  55  (2000),  no. 3, 475--546.

%\bibitem{It} Itenberg, I.,
%Amibes de vari\'et\'es alg\'ebriques et d\'enombrement de courbes [d'apr\`es G. Mikhalkin],
%In {\em S\'eminaire Bourbaki (2002/03)},
%Ast\'erisque  No. 294 (2004), 335--361.

\bibitem{Ka} Kapranov, M.,
Amoebas over non-Archimedian fields, Preprint 2000.

\bibitem{IKS}  Itenberg, I.,
Kharlamov, V., Shustin, E.,
Welschinger invariant and enumeration of real rational curves.
{\em Int. Math. Res. Not.}  2003,  no. 49, 2639--2653.

\bibitem{KM}  Kontsevich, M., Manin, Yu.,
{Gromov-Witten classes, quantum cohomology and
enumerative geometry}. {\em Comm. Math. Phys.} {\bf 164} (1994), 525-562.

\bibitem{KS} Kontsevich, M., Soibelman, Ya.,
Homological mirror symmetry and torus fibrations.
In {\em Symplectic geometry and mirror symmetry (Seoul, 2000)},
World Sci. Publishing, River Edge, NJ, 2001, 203--263.

\bibitem{KS-K3} Kontsevich, M., Soibelman, Ya.,
Affine structures and non-archimedean analytic spaces,
In {\em The unity of mathematics} Progr. Math., {\bf 244},
Birkh\"auser Boston, Boston, MA, 2006,  321--385.

\bibitem{LiMa}  Litvinov, G. L.,
The Maslov dequantization, idempotent and tropical mathematics: a very brief introduction.
In {\em Idempotent mathematics and mathematical physics},
Contemp. Math., 377, Amer. Math. Soc., Providence, RI, 2005, 1--17.

\bibitem{Mi-am} Mikhalkin, G.,
Real algebraic curves, the moment map and amoebas.
{\em Ann. of Math. (2)},  151  (2000),  no. 1, 309--326.

\bibitem{Mi-pp} Mikhalkin, G.,
Decomposition into pairs-of-pants for complex algebraic hypersurfaces.
{\em Topology}  43  (2004),  no. 5, 1035--1065.

\bibitem{Mi-GV}  Mikhalkin, G.,
Amoebas of algebraic varieties and tropical geometry.
In {\em Different faces of geometry,}
Int. Math. Ser. (N. Y.), Kluwer/Plenum, New York, 2004, 257--300.

\bibitem{Mi-R2}  Mikhalkin, G.,
Enumerative tropical algebraic geometry in $\Bbb R\sp 2$.
{\em J. Amer. Math. Soc.}  18  (2005),  no. 2, 313--377.

\bibitem{Mi-Q} Mikhalkin G., Rational tropical curves in $\R^n$,
to appear.

\bibitem{M-Zeut} Mikhalkin G., Zeuthen's numbers for toric surfaces
via tropical geometry, to appear.

\bibitem{MiRu} Mikhalkin, G., Rullgard, H.,
Amoebas of maximal area.
{\em Internat. Math. Res. Notices},  2001,  no. 9, 441--451.

\bibitem{MZ} Mikhalkin, G., Zharkov, I.,
Tropical curves, their Jacobians and Theta-functions, to appear.

\bibitem{NS} Nishinou, T., Siebert, B.,
{\em Toric degenerations of toric varieties and tropical curves},
http://arxiv.org/abs/math.AG/0409060.

\bibitem{RGST} Richter-Gebert, J., Sturmfels, B., Theobald, Th.,
First steps in tropical geometry.
In {\em Idempotent mathematics and mathematical physics,}
Contemp. Math., 377, Amer. Math. Soc., Providence, RI, 2005, 289--317.

%\bibitem{Vi} Risler, J.-J., Construction d'hypersurfaces r\'eelles [d'apr\`es Viro].
%In  {\em S\'eminaire Bourbaki (1992/93)}, Aste'risque, vol. 216, p. 69--86.

\bibitem{PR}  Passare, M., Rullgard, H.,
Amoebas, Monge-Amp\`ere measures, and triangulations of the Newton polytope.
{\em Duke Math. J.}  121  (2004),  no. 3, 481--507.

\bibitem{PST} Passare, M., Sadykov, T., Tsikh, A.,
Singularities of hypergeometric functions in several variables.
{\em Compos. Math.}  141  (2005),  no. 3, 787--810.

\bibitem{Sh} Shustin, E.,
Patchworking singular algebraic curves, non-Archimedean amoebas and enumerative geometry,
http://arxiv.org/abs/math.AG/0211278.
%A tropical approach to enumerative geometry.
%{\em Algebra i Analiz}  17  (2005),  no. 2, 170--214.

\bibitem{Sp-diss} Speyer, D.,
Dissertation, University of California, Berkeley, 2005.

\bibitem{SS} Speyer, D., Sturmfels, B.,
The tropical Grassmannian.  {\em Adv. Geom.}  4  (2004),  no. 3, 389--411.

\bibitem{S9} Sturmfels, B.,
{\em Solving systems of polynomial equations.}
CBMS Regional Conference Series in Mathematics, 97,
American Mathematical Society, Providence, RI, 2002.

\bibitem{Vi} Viro, O. Ya., Gluing of algebraic hypersurfaces,
smoothing of singularities and construction of curves.
In {\em Proc. Leningrad Int. Topological Conf., Leningrad, 1982},
Nauka, Leningrad, 1983, 149--197.

\bibitem{Vi-deq}  Viro, O. Ya., Dequantization of real algebraic geometry on logarithmic paper.
In {\em European Congress of Mathematics, Vol. I (Barcelona, 2000),}
Progr. Math., 201, Birkh\"auser, Basel, 2001, 135--146.

\bibitem{We} Welschinger, J.-Y., {Invariants of real rational
symplectic 4-manifolds and lower bounds in real enumerative geometry.}
{\em C. R. Math. Acad. Sci. Paris}  {\bf 336}  (2003),  no. 4, 341--344.

\bibitem{We3} Welschinger, J.-Y.,
Spinor states of real rational curves in real algebraic convex 3-manifolds and enumerative invariants,
http://arxiv.org/abs/math.AG/0311466.

\end{thebibliography}
\end{document}